\documentclass[a4paper,11pt]{article}
\usepackage{amssymb,palatino}
\newtheorem{theorem}{Theorem}
\newtheorem{corollary}{Corollary}
\newtheorem{lemma}{Lemma}

\newtheorem{proposition}{Proposition}
\newcommand{\sty}{\displaystyle}
\newcommand{\QED}{\begin{flushright} $\Box$ \end{flushright}}
\newcommand{\Tr}{{\rm Tr} }                % traccia
\newcommand{\jacobi}[2]{\left( \begin{array}{c} #1\\ \hline #2 \end{array} \right)}
\begin{document}

\title{Gauss Sums of Cubic Characters over $\mathbb F_{p^r}$, $p$ odd}% prime}

\author{Michele Elia \thanks{Politecnico di Torino, Italy}, ~~
 Davide Schipani\thanks{University of Zurich, Switzerland}
}

%\date{August 2010}
\maketitle

\thispagestyle{empty}

\begin{abstract}
\noindent
An elementary approach is shown which derives the values of the Gauss sums over $\mathbb F_{p^r}$, $p$ odd, 
 of a cubic character. New links between
 Gauss sums over different field extensions are shown in terms of factorizations of the Gauss sums themselves, which then are revisited in terms of
 prime ideal decompositions.
%\noindent
Interestingly, one of these results gives a representation of primes $p$ of the form $6k+1$ by
 a binary quadratic form in integers of a subfield of the cyclotomic field of
 the $p$-th roots of unity.  
\end{abstract}

\vspace{2mm}

% ********************************************************************
\vspace{8mm}

\section{Introduction}
Let $\mathbb F_{p^{r}}$ be a Galois field of order $p^r$,
 with $\Tr_r(x)=\sum_{j=0}^{r-1} x^{p^j}$
 being the trace function over $\mathbb F_{p^r}$, and $\Tr_{r/d}(x)=\sum_{j=0}^{r/d-1} x^{p^{dj}}$
 the relative trace function over $\mathbb F_{p^r}$ relatively to $\mathbb F_{p^d}$, with $d|r$
 \cite{lidl}. 

\noindent
Further let $\chi_m$ be a character of order $m$ defined over $\mathbb F_{p^r}$ and
 taking values in the cyclotomic field  $\mathbb Q(\zeta_m)$, where $\zeta_m$ denotes
 a primitive $m$-th complex root of unity.

\noindent 
The Gauss sum of $\chi_m$ over $\mathbb F_{p^{r}}$ is defined, \cite{berndt,rosen}, for any $\beta\in\mathbb F_{p^r}$ as 
$$  G_r(\beta, \chi_m) = \sum_{y \in  \mathbb F_{p^r}}  \chi_m(y) \zeta_p^{\Tr_r(\beta y)} = 
    \bar \chi_m(\beta) G_r(1, \chi_m) ~~. $$
We will focus our interest on cubic characters $\chi_3$ for odd primes $p$, while the case $p=2$ is dealt with in \cite{gauss1}. A cubic character $\chi_3$ can be either the principal character, i.e. 
 $\chi_3(\beta)=1$ for all $\beta \in \mathbb F_{p^{r}}^*$, or a non-principal character (if $p^r \equiv 1 \bmod 6$)
 
 $$  \chi_3(\alpha^{h+3j}) = \zeta_3^{h} ~~~~h=0,1,2~,~j\in\mathbb{N}~~,  $$ 
 where $\alpha$ is a generator of $\mathbb F_{p^{r}}^*$. 
 In addition, $\chi_3(0)=0$ by definition. \\
By the above assumptions, the values of the Gauss sums of a cubic character over $\mathbb F_{p^{r}}$
 are in general algebraic integers in the field $\mathbb Q(\zeta_3,\zeta_p)=\mathbb Q(\zeta_{3p})$,
 $p \neq 3$ and $\mathbb Q(\zeta_3)$, if $p=3$. 
Our aim is to give a thorough overview and derive more precise statements about these values with elementary techniques.
In particular, we have obtained an interesting new result in the vein of Gauss' closed expression for $G_1(1, \chi_2)$ in terms of a fourth root of unity and a root of an integral quadratic polynomial.
% in terms of $i=\zeta_4$ and a root of $x^2-p$. 
%
Specifically, our equation (\ref{gauss3}) expresses 
  $G_1(1, \chi_3)$ in terms of a cubic root of unity and a root of an
  integral cubic polynomial (see Equation (\ref{gausscubic})),   
  whose coefficients are given explicit functions of $p$.
%  , which can be computed in polynomial complexity.
 Gauss sums over extended fields are obtained from the expression of Gauss sums over %the lowest field 
 %(possibly $\mathbb F_p$) of a field tower is known, the  can be obtained 
smaller fields, either using Davenport-Hasse's theorem or with new methods developed here to link values over different field extensions.

%tegers of a subfield of the cyclotomic field of
 %the $p$-th roots of unity.  
  
%Namely, we derive expressions for the Gauss sums over $\mathbb F_{p}$, which are then used to obtain the %Gauss sums over extended fields (either by means of Davenport-Hasse's theorem or with new techniques %developed here).
 %(alternatively to Davenport-Hasse's theorem).

\section{Lemmas}
For the considerations below, let
$$ A_{s}(\alpha) = \sum_{y \in  \mathbb F_{p^{s}}}  \chi_m(y+\alpha) ~~~~\mbox{and}~~~~ 
B_{r,s}(\alpha) =\sum_{z_1 \ldots , z_{r-1} \in  \mathbb F_{p^{s}}} \chi_m(1+ \sum_{i=1}^{r-1} z_i\alpha^i)   ~~, $$
where $\alpha$ may be in an extension of $\mathbb F_{p^{s}}$, and $\chi_m$ defined in the same extension.

\begin{lemma}
  \label{lemma1pbis}
Let $\chi_m$ be a nontrivial character of order $m$ over $\mathbb F_{p^{rs}}$, $p$ prime, whose restriction to $\mathbb F_{p^{s}}$ is also nontrivial, and assume %that $r$ divides $p^s-1$; 
 that there exists an irreducible polynomial $X^r-\beta$, for a suitable $\beta \in  \mathbb F_{p^s}$. %(that is $\chi_{\Lambda(r)}(\beta)=\zeta_{\Lambda(r)}$ for a character of order $r$ over $\mathbb F_{p^s}$)
Then
$$  G_{rs}(1, \chi_m)= \bar \chi_m(r)G_s(1, \chi_m) B_{r,s}(\alpha) , $$
where $\alpha$ is a root of $X^r-\beta$ (thus with relative trace $\Tr_{rs/s}(\alpha)=0$).
\end{lemma}

\noindent
{\sc Proof}.
Since $\mathbb F_{p^{rs}}$ is an extension of order $r$ of $\mathbb F_{p^{s}}$, its elements can be written in the form $\sum_{i=0}^{r-1} x_i\alpha^i$ with $x_0, \ldots , x_{r-1} \in \mathbb F_{p^{s}}$. We thus have 
$$ \begin{array}{lcl}
    G_{rs}(1, \chi_m) &=& \sty \sum_{z \in  \mathbb F_{p^{rs}}}   \chi_m(z) \zeta_p^{\Tr_{rs}(z)} = 
   \sum_{x_0, \ldots , x_{r-1}  \in  \mathbb F_{p^{s}}}   \chi_m(\sum_{i=0}^{r-1} x_i\alpha^i) \zeta_p^{\Tr_{rs}(\sum_{i=0}^{r-1} x_i\alpha^i)} \\
    &=& \sty \sum_{x_0, \ldots , x_{r-1} \in  \mathbb F_{p^{s}}} \chi_m(\sum_{i=0}^{r-1} x_i\alpha^i) \zeta_p^{\Tr_{s}(rx_0)} =  \sty \sum_{\stackrel{x_0 \in F_{p^{s}}^*}
     {x_1 \ldots , x_{r-1} \in  \mathbb F_{p^{s}}}} \chi_m(\sum_{i=0}^{r-1} x_i\alpha^i) \zeta_p^{\Tr_{s}(rx_0)},
    \end{array}   ~~ $$ 
where we used that
$$\Tr_{rs}(\sum_{i=0}^{r-1} x_i\alpha^i)=\Tr_{s}(\Tr_{rs/s}(\sum_{i=0}^{r-1} x_i\alpha^i)) =
      \sum_{i=0}^{r-1} \Tr_{s}(x_i\Tr_{rs/s}(\alpha^i)) = \Tr_{s}(rx_0) ~,  $$
since, if $\alpha$ is a solution of $X^r-\beta=0$, then $\Tr_{rs/s}(\alpha^j)=0$ for every $j=1,\ldots, r-1$ for example as a consequence of the Newton formulas \cite[vol. I, pg. 166]{burnside}. 
%
%as $\Tr_{rs/s}(x_i\alpha^i)=x_i\Tr_{rs/s}(\alpha^i)$ because $x_i^{p^s}=x_i$, further $\Tr_{rs/s}(1)=r$ 
 %and $\Tr_{rs/s}(\alpha^i)=0$ if $i >0$ because $\alpha$ is a root of the polynomial $X^r-\beta$ which
 % has every coefficient (that is an elementary symmetric function of the roots) equal to zero with the exclusion of the constant coefficient, then every symmetric function of the powers of the roots
 % $s_i= \sum_{j=0}^{r-1} (\alpha^{p^{js}})^i= \Tr_{rs/s}(\alpha^i)$ with $i=1, \ldots, r-1$
 % is zero as a consequence of the Newton formulas \cite[vol. I, pg. 166]{burnside},
 % we thus have $\Tr_{rs}(\sum_{i=0}^{r-1} x_i\alpha^i) =  \Tr_{s}(rx_0)$;
The sum  $\sum_{x_1 \ldots , x_{r-1} \in  \mathbb F_{p^{s}}} \chi_m(\sum_{i=1}^{r-1} x_i\alpha^i) $ is zero, since with the change of variable $x_i=x_i'\lambda$, where $\lambda$ is an element of $\mathbb F_{p^{s}}^*$ with $\chi_m(\lambda)\neq 1$, the sum becomes $\chi_m(\lambda)\sum_{x_1 \ldots , x_{r-1} \in  \mathbb F_{p^{s}}} \chi_m(\sum_{i=1}^{r-1} x_i\alpha^i) $. \\
Now, as $x_0 \neq 0$, we may perform the change of variables $x_i = x_0 z_i$, $i=1, \ldots , r-1$ and write 
$$ \begin{array}{lcl}
 G_{rs}(1, \chi_m) &=& \sty\sum_{x_0 \in F_{p^{s}}^*} \zeta_p^{\Tr_{s}(rx_0)} \sum_{z_1 \ldots , z_{r-1} \in  \mathbb F_{p^{s}}} \chi_m(x_0 + x_0 \sum_{i=1}^{r-1} z_i\alpha^i) \\
  & =& \sty \sum_{x_0 \in F_{p^{s}}^*} \zeta_p^{\Tr_{s}(rx_0)} \chi_m(x_0)\sum_{z_1 \ldots , z_{r-1} \in  \mathbb F_{p^{s}}} \chi_m(1+ \sum_{i=1}^{r-1} z_i\alpha^i) ~~.
    \end{array}   ~~ $$ 
The conclusion is immediate, noting that the first summation is simply $\bar \chi_m(r) G_s(1, \chi_m)$.  \hfill $\Box$
%\QED

\vspace{5mm}
\noindent
The following Lemma is a corollary of the previous one, specialized to the case $r=2$.
 However, we present another proof, whose structure and running results are instrumental
 to proofs of further theorems. 

\begin{lemma}
  \label{lemma1p}
Let $\chi_m$ be a nontrivial character of order $m$ over $\mathbb F_{p^{2s}}$, $p$ odd, whose restriction to $\mathbb F_{p^{s}}$ is also nontrivial; %, with $p=6k+5$;
 then
$$  G_{2s}(1, \chi_m)= \chi_m(\alpha)G_s(1, \chi_m) A_s(\frac{1}{2\alpha}) ~~, $$
where $\alpha$ is a root of an irreducible polynomial $X^2-\beta$ for a suitable $\beta \in  \mathbb F_{p^s}$. %(thus with relative trace $\Tr_{2s/s}(\alpha)=0$),  (that is $\chi_2(\beta)=-1$ for a quadratic character over $\mathbb F_{p^s}$).
\end{lemma}
We note that from the definition of $\alpha$ and $\beta$ it follows that $\Tr_{2s/s}(\alpha)=0$ and that $\chi_2(\beta)=-1$ for the nontrivial quadratic character over $\mathbb F_{p^s}$. 

\noindent
{\sc Proof}.
%We note that from the definition of $\alpha$ and $\beta$ it follows that $\Tr_{2s/s}(\alpha)=0$ and that $\chi_2(\beta)=-1$ for the nontrivial quadratic character over $\mathbb F_{p^s}$. \\
%
Since $\mathbb F_{p^{2s}}$ is a quadratic extension of $\mathbb F_{p^{s}}$, its elements can be written in the form $x+\alpha y$ with $x,y, \in \mathbb F_{p^{s}}$. We thus have 
$$ G_{2s}(1, \chi_m)= \sum_{z \in  \mathbb F_{p^{2s}}}   \chi_m(z) \zeta_p^{\Tr_{2s}(z)} = 
   \sum_{x,y \in  \mathbb F_{p^{s}}}   \chi_m(x+\alpha y) \zeta_p^{\Tr_{2s}(x+\alpha y)}
   =  \sum_{x,y \in  \mathbb F_{p^{s}}}   \chi_m(x+\alpha y) \zeta_p^{\Tr_{s}(2x)} ~~, $$ 
where we have used the equality $\Tr_{2s}(x+\alpha y) = 2 \Tr_{s}(x)= \Tr_{s}(2x)$.
%since 
%$$\Tr_{2s}(x+\alpha y)=\Tr_{2s}(x)+\Tr_{2s}(\alpha y) =
%      2 \Tr_{s}(x)+ \Tr_{s}(\alpha y) + \Tr_{s}(\alpha^{p^s} y) = 2 \Tr_{s}(x)+ \Tr_{s}(\alpha y) + \Tr_{s}(-\alpha y) ~,  $$
%as the Frobenius map $x\to x^{p^{s}}$ permutes the roots of $X^2-\beta$. We thus have
% $\Tr_{2s}(x+\alpha y) = \Tr_{2s}(x)= \Tr_{s}(2x)$.
 Multiplying the last sum by $\bar \chi_m(2) \chi_m(2) =1$, we can write
$$ G_{2s}(1, \chi_m)= \bar \chi_m(2) \sum_{x',y \in  \mathbb F_{p^{s}}}   \chi_m(x'+2\alpha y) \zeta_p^{\Tr_{s}(x')} ~, $$
and split the summation into three sums
$$ \bar \chi_m(2) \sum_{y \in  \mathbb F_{p^{s}}}   \chi_m(2\alpha y) ~~,~~
 \bar \chi_m(2) \sum_{x' \in  \mathbb F_{p^{s}}^*}   \chi_m(x') \zeta_p^{\Tr_{s}(x')} ~~,~~
\bar \chi_m(2) \sum_{x',y \in  \mathbb F_{p^{s}}^*}   \chi_m(x'+2\alpha y) \zeta_p^{\Tr_{s}(x')} ~~.$$
The first summation is $0$, the second summation is $\bar \chi_m(2) G_s(1, \chi_m)$; the third summation
 can ~ be written as follows: the substitution $y=z x'$ yields
$$ \bar \chi_m(2) \sum_{x',z \in  \mathbb F_{p^{s}}^*}   \chi_m(x'+2\alpha z x') \zeta_p^{\Tr_{s}(x')}=
    \bar \chi_m(2) \sum_{x' \in  \mathbb F_{p^{s}}^*}   \chi_m(x') \zeta_p^{\Tr_{s}(x')} 
                 \sum_{z \in  \mathbb F_{p^{s}}^*}  \chi_m(1+2\alpha z) = $$
$$     \bar \chi_m(2) G_s(1, \chi_m) \chi_m(2\alpha) \sum_{z \in  \mathbb F_{p^{s}}^*}  \chi_m(z+\frac{1}{2\alpha}) =
     \bar \chi_m(2) G_s(1, \chi_m) \chi_m(2\alpha) [A_s(\frac{1}{2\alpha})- \chi_m(\frac{1}{2\alpha})] ~~.$$
In conclusion, by combining the above summations, we have
$ G_{2s}(1, \chi_m)= \chi_m(\alpha)G_s(1, \chi_m) A_s(\frac{1}{2\alpha})$. 
\QED

\begin{corollary}
Suppose $p$ is odd and $t=2^ks$, with $k\geq 1$, and let $\chi_m$ be a nontrivial character over $\mathbb F_{p^{t}}$, whose restriction to $\mathbb F_{p^{s}}$ is also nontrivial. Then 
$$
G_t(1,\chi_m)=G_{s}(1,\chi_m)\prod_{i=1}^{k} \chi_m(\alpha_i)A_{2^{i-1}s}(\frac{1}{2\alpha_i}),
$$
where $\alpha_i$ is a root of an irreducible polynomial $X^2-\beta_i$ over $\mathbb F_{p^{2^{i-1}s}}$, $i=1,\ldots,k$.
\end{corollary}

\begin{lemma}
  \label{lemma2}
Let $\chi_m$ be a character over $\mathbb F_{p^{s}}$,
 $p\equiv -1$ mod $m$ and $m$ odd. Then $G_{s}(1, \chi_m)$ is real.
\end{lemma}

\noindent
{\sc Proof}. 
We can write 
$$ G_{s}(1, \chi_m)= G_0 +\zeta_m G_1 + \zeta_m^2 G_2+\cdots +\zeta_{m}^{m-1} G_{m-1} ~~,$$
 where $\zeta_m$ is a primitive $m$-th root of unity and
$$G_j= \sum_{\chi_m( x)=\zeta_m^j} \zeta_p^{\Tr_{s}( x)} ~~,~~0\leq j\leq m-1~~, $$
known as Gauss periods \cite{gauss}, are real numbers, since $\chi_m(x)=\chi_m(-x)$, as $-1=(-1)^m$ is an $m$-th power. %given that $$-1=-\gamma^{p^{s}-1}=(-\gamma^{\frac{p^{s}-1}{m}})^m $$
 %with $\gamma$ a primitive element of $\mathbb F_{p^{s}}$
Thus, in each sum the exponentials  occur in complex conjugated pairs.
Furthermore, $G_j=G_{m-j}$ as proved by the following chain of equalities:
$$
G_j=\sum_{\chi_m(x)=\zeta_m^j} \zeta_p^{\Tr_{s}(x)}=\sum_{\chi_m(x)=\zeta_m^j}
     \zeta_p^{\Tr_{s}( x^p)}=\sum_{\chi_m(x^p)=\zeta_m^{pj}} \zeta_p^{\Tr_{s}( x^p)}=
     \sum_{\chi_m( y)=\zeta_m^{m-j}} \zeta_p^{\Tr_{s}(y)}=G_{m-j} ~~.
$$
In fact raising the trace argument to the power $p$ leaves the trace invariant; $\zeta_m^{pj}=\zeta_m^{-j}$ as $p$ is congruent to $-1$ modulo $m$; lastly, the automorphism $\sigma(x)= x^p$ simply permutes the elements of the field. 
Then, for any $j$, $\zeta_{m}^{j} G_{j}$ and $\zeta_{m}^{m-j} G_{m-j}$ sum to give a real number, %(if
 %$m$ is even and $j=m-j$ then $\zeta_m^j$ is already real)
 hence $ G_{s}(1, \chi_m)$ is also real.
\QED 

\begin{corollary}
Let $\chi_3$ be a cubic character, $p\equiv 2$ modulo $3$ ($p=2$ or $p=6k+5$).
 Then $G_{s}(1, \chi_3)$ is real.
\end{corollary}

\paragraph{Remark 1.}
For the case $p=2$, see an alternative proof in \cite{gauss1}.

\begin{lemma}
  \label{lemma3}
%Let $\chi_m$ be a non-trivial character over $\mathbb F_{p^{2s}}$,
% $p\equiv -1$ mod $m$ with $m$ and $s$ odd. If the restriction of $\chi_m$ 
 %to $\mathbb F_{p}$ is trivial, then $G_{2s}(1, \chi_m)=p^s$.
Let $\chi_m$ be a nontrivial character over $\mathbb F_{p^{2s}}$,
 with $p$ and $m$ odd. If $(p^s-1,m)=1$ (in particular the restriction of $\chi_m$ 
 to $\mathbb F_{p^s}$ is trivial), then $G_{2s}(1, \chi_m)=p^s$.
\end{lemma}

\noindent
{\sc Proof}. 
As in Lemma \ref{lemma1p}, let $\alpha$ be defined as a root of an irreducible polynomial $X^2-\beta$, with a suitable $\beta \in  \mathbb F_{p^{s}}$. Then %the Gauss sum $G_{2s}(1, \chi_m)$
 %can be written as
$$ G_{2s}(1, \chi_m)= %\sum_{z \in  \mathbb F_{p^{2s}}}   \chi_m(z) \zeta_p^{\Tr_{2s}(z)} = 
   %\sum_{x,y \in  \mathbb F_{p^{s}}}   \chi_m(x+\alpha y) \zeta_p^{\Tr_{2s}(x+\alpha y)}=
    \sum_{x,y \in  \mathbb F_{p^{s}}}   \chi_m(x+\alpha y) \zeta_p^{\Tr_{s}(2x)} ~~. $$  
We split the summation into three: $S_1=\bar \chi_m(2) \sum_{y \in  \mathbb F_{p^{s}}}   \chi_m(2\alpha y)$,
$$  S_2=\bar \chi_m(2) \sum_{x' \in  \mathbb F_{p^{s}}^*}   \chi_m(x') \zeta_p^{\Tr_{s}(x')}, ~~\mbox{and}~~
S_3=\bar \chi_m(2) \sum_{x',y \in  \mathbb F_{p^{s}}^*}   \chi_m(x'+2\alpha y) \zeta_p^{\Tr_{s}(x')} ~~.$$
The first summation is $S_1=\chi_m(\alpha)(p^s-1)$, since the character is trivial over $\mathbb F_{p^{s}}$, the second summation is $S_2=-\bar \chi_m(2)$, and the third summation, after the substitution $y=z x'$,  gives
$$ S_3=
%\bar \chi(2) \sum_{x',z \in  \mathbb F_{p^{s}}^*}   \chi(x'+2\alpha z x') e^{\frac{2\pi i}{p} \Tr_{s}(x')}=
    \bar \chi_m(2) \sum_{x' \in  \mathbb F_{p^{s}}^*}   \chi_m(x') \zeta_p^{\Tr_{s}(x')} 
      \sum_{z \in  \mathbb F_{p^{s}}^*}  \chi_m(1+2\alpha z) = % $$ $$
       -\bar \chi_m(2) [  \chi_m(2\alpha) \sum_{z \in  \mathbb F_{p^{s}}}  \chi_m(\frac{1}{2\alpha}+z) - 1] ~~.$$
In order to evaluate $A_s(\frac{1}{2\alpha})= \sum_{z \in  \mathbb F_{p^{s}}}  \chi_m(\frac{1}{2\alpha}+z)$, we consider the sum of $A_s(\beta)$, for every $\beta \in \mathbb F_{p^{2s}}$, and observe that $A_s(\beta)=p^s-1$ if $\beta \in \mathbb F_{p^s}$, since all elements in this field are $m$-th powers, while,
 if $\beta \not \in \mathbb F_{p^s}$ all sums assume the same value $A_s(\alpha)$, which is shown as follows: set $\beta=u+\alpha v$ with $v \neq 0$, then
$$ \sum_{z \in  \mathbb F_{p^s}} \chi_m(z+u+ \alpha v )  = \sum_{z \in  \mathbb F_{p^s}} \chi_m(v) \chi_m((z+u)v^{-1}+ \alpha )
    = \sum_{z' \in  \mathbb F_{p^s}} \chi_m(z'+ \alpha )=A_s(\alpha)~~. $$
Therefore, the sum  $\displaystyle \sum_{\beta \in  \mathbb F_{p^{2s}}} A(\beta) = \sum_{\beta \in  \mathbb F_{p^{2s}}}\sum_{z \in  \mathbb F_{p^s}} \chi_m(z+\beta)=\sum_{z \in  \mathbb F_{p^s}}\sum_{\beta \in  \mathbb F_{p^{2s}}} \chi_m(z+\beta)=0$ yields
$$  p^s(p^s-1)+(p^{2s}-p^s) A(\alpha) =0 $$
which implies $A(\alpha) =-1=A(\frac{1}{2\alpha})$. %(clearly $\frac{1}{2\alpha}$ is not in $\mathbb F_{p^s}$, otherwise $\alpha=\frac{1}{2\alpha}2\alpha^2=2\beta\frac{1}{2\alpha}$ would also be in $\mathbb F_{p^s}$).
Finally, by combining the above,
$$ G_{2s}(1,\chi_m)= \chi_m(\alpha) (p^s-1) + \chi_m(\alpha) = \chi_m(\alpha) p^s= p^s~~, $$
because  $\alpha$, a root of $X^2-\beta$, is an $m$-th power, since every $\beta \in \mathbb F_{p^s}$ is an $m$-th power.
 %\hfill $\Box$
\QED

\paragraph{Remark 2.}
The above lemma can also be proved using a theorem by Stickelberger,
 \cite[Theorem 5.16]{lidl} or \cite{stickelberger}.

\section{Results}
%The case of trivial characters is reported here for the sake of completeness.

\paragraph{Trivial character.} Let $\chi_3$ be trivial, then
$$ G_r(1, \chi_3) = \sum_{y \in  \mathbb F_{p^r}}  \chi_3(y) \zeta_p^{\Tr_r( y)} = 
  \sum_{y \in  \mathbb F_{p^r}} \zeta_p^{\Tr_r( y)}-1 = 
  p^{r-1} \sum_{a \in  \mathbb F_{p}} \zeta_p^{a}-1 = -1 ~~, $$
since the number of elements with the same trace $a \in \mathbb F_p$
 ($0$ included) is equal to $p^{r-1}$, i.e. the number of roots in $\mathbb F_{p^r}$ of the equation
  $\Tr_r(x)=a$.  This result settles in particular all the cases of the fields $ \mathbb F_{3^r}$, or 
  $ \mathbb F_{p^r}$ with $p\equiv 5 \bmod 6$ and odd $r$, where there is only 
  the principal character, because every field element is a cube. 
%$\beta$ is a cube. %as the following chain of equalities shows
%$$ \beta \cdot 1 =  \beta \cdot (\beta^{p^r-1})^2= \beta \beta^{2p^{r}-2} = \beta^{2p^{r}-1} =(\beta^{\frac{2p^{r}-1}{3}})^3 ~~,  $$   
%since $\beta^{p^r-1}=1$, and $p = -1 \bmod 3$, so that $2p^{r}-1$ is divisible by $3$.

%If $p=3$, then, for any $r$, $3^r-1$ is not divisible by $3$ and only the trivial character exists. 
%Every $\beta$ in $\mathbb F_{3^r}$ is a cube, as 
%$$ \beta \cdot 1 =  \beta \cdot (\beta^{3^r-1})= (\beta^{3^{r-1}})^{3} ~~,  $$   
%since $\beta^{3^r-1}=1$. \\
%
%In this case, since $ \mathbb F_{3^r}$ has only the trivial cubic character $\chi_3(x)=1$, $x\neq 0$,
% we have
%$$  G_1(1,\chi_3)  = \sum_{y \in  \mathbb F_{3}^*}  \zeta_3^{y} = \sum_{y \in  \mathbb F_{3}}  \zeta_3^{y} -1 = -1 ~~, $$
%  as $\sum_{a \in  \mathbb F_{3}} \zeta_3^{a}=0$. This same property
 %holds for every $r > 1$; in general we have
%$$ G_r(1, \chi_3) = \sum_{y \in  \mathbb F_{3^r}}  \chi_3(y) \zeta_3^{\Tr_r( y)} = 
%  \sum_{y \in  \mathbb F_{3^r}} \zeta_3^{\Tr_r( y)}-1 = 
%  3^{r-1} \sum_{a \in  \mathbb F_{3}} \zeta_3^{a}-1 = -1 ~~, $$
%as the character $\chi_3$ is trivial and the number of elements with the same trace $a \in \mathbb F_3$
 %($0$ included) is equal to $3^{r-1}$, i.e. the number of roots in $\mathbb F_{3^r}$ of the equation
%  $\Tr_r(x)=a$.  

\paragraph{Nontrivial character: case p=6k+5.}

If $p=6k+5$ and $r$ is even, a nontrivial cubic character exists and it will be shown that
 $G_r(1, \chi_3) = -(-p)^{r/2}$, without recurring to Davenport-Hasse's theorem.

\begin{theorem}
  \label{lemma2p}
If $p=6k+5$ and $s$ is odd, then $G_{2s}(1, \chi_3)=p^s$.
\end{theorem}

\noindent
{\sc Proof}.
Since $p^s=-1 \bmod 3$, the conclusion is a consequence of Lemma \ref{lemma3}.
\QED

\begin{theorem}
  \label{lemma3p}
If $p=6k+5$ and $s$ is even, then $G_{2s}(1, \chi_3)=(-p)^{s/2}G_{s}(1, \chi_3)$.
\end{theorem}

\noindent
{\sc Proof}.
Let $\alpha \in \mathbb F_{p^{2s}}$ be a cube and root of an irreducible polynomial $X^2-\beta$ over
 $\mathbb F_{p^{s}}$ (clearly such an $\alpha$ exists, since if $\gamma$ is a root of 
 $X^2-\beta$, with $\chi_2(\beta)=-1$, then $\gamma^3$, a cube, is a root of $X^2-\beta^3$ and
 $\chi_2(\beta^3)=\chi_2(\beta)^3=-1$). 
Then by Lemma \ref{lemma1p}
$$  G_{2s}(1, \chi_3)= G_{s}(1, \chi_3) A_{s}(\frac{1}{2\alpha})~~,  $$
where $A_{s}(\frac{1}{2\alpha})= \sum_{z \in \mathbb F_{p^{s}}} \chi_3(\frac{1}{2\alpha}+z)$ is
 an algebraic integer in the cyclotomic field $\mathbb Q(\zeta_3)$ which can be written as
 $A_0+\zeta_3 A_1 + \zeta_3^2 A_2$, where $A_0$, $A_1$, and $A_2$ are the numbers of $z$ for which  $\chi_3(\frac{1}{2\alpha}+z)$ is equal to
 $1$, $\zeta_3$ or $\zeta_3^2$, respectively, and $A_0+A_1+A_2=p^s$. 

Now, by Lemma \ref{lemma2}, both $ G_{2s}(1, \chi_3)$ and  $G_{s}(1, \chi_3)$ are real, which implies that $A_{s}(\frac{1}{2\alpha})$ is also real, so that $A_1=A_2$. We also know that $A_0+A_1+A_2=A_0+2A_1=p^s$,
 so we consider two equations for $A_0$ and $A_1$:
$$   \left\{ \begin{array}{l}
       A_0 + 2 A_1 = p^s  \\
       A_0-A_1 = \pm p^{s/2}  \\
             \end{array}   \right.
$$
obtained from the fact that we know the absolute values of $ G_{s}(1, \chi_3)$ and $ G_{2s}(1, \chi_3)$, \cite[Theorem 1.1.4, pg. 10]{berndt}, \cite{gauss1}.\\
Solving for $A_1$ we have $A_1= \frac{1}{3}(p^s \mp  p^{s/2}) $. As $A_1$ must be an integer, we have
$$  A_{s}(\frac{1}{2\alpha})= A_0 - A_1 =  \left\{ \begin{array}{lcl}
        p^{s/2}    &~~& \mbox{if $s/2$ is even}\\
        -p^{s/2}&~~& \mbox{if $s/2$ is odd}.  \\
             \end{array}   \right.  ~~
$$
\QED

\begin{corollary}
If $p=6k+5$ and $s$ is even, then $G_{2s}(1, \chi_3)=-p^s$.
\end{corollary}

\paragraph{Nontrivial character: case p=6k+1.}
If $p=6k+1$, $p^r-1$ is divisible by $3$, so there exists a nontrivial cubic character in $\mathbb F_{p^r}$
for every $r \geq 1$: 
we know that the Gauss sum over $\mathbb F_{p}$ of a nontrivial cubic character is an algebraic integer in
 $\mathbb Q(\zeta_{3p})$ of absolute value $\sqrt p$.
 Specifically we have
\begin{theorem}
   \label{theo3} 
If $p=6k+1$, then $G_{1}(1, \chi_3)$ is an element of $\mathbb Q(\zeta_{3}, \eta)$, a subfield of
 $\mathbb Q(\zeta_{3p})$ with degree $6$ over $\mathbb Q$, where $\eta$ is a root of a cubic 
 polynomial with rational integer coefficients and cyclic Galois group over $\mathbb Q$. 
\end{theorem} 

\noindent
{\sc Proof}.
As in the proof of Lemma \ref{lemma2}, for a cubic character we can write
$$ G_{1}(1, \chi_3)= G_0 +\zeta_3 G_1 + \zeta_3^2 G_2 ~~,$$
 where $\zeta_3$ is a primitive cube root of unity and, for $0\leq j\leq 2$,
 $$G_j= \sum_{\chi_3(x)=\zeta_3^j} \zeta_p^{x}~~, $$
which are real numbers since $\chi_3(x)=\chi_3(-x)$, as $-1$ is a cubic power.
Then, to evaluate the Gauss sum $G_{1}(1, \chi_3)$ is tantamount to computing the Gauss periods
$G_0$, $G_1$, and $G_2$.
 The following derivation can be found partly, in different 
 form, in Gauss \cite[art. 350-352]{gauss}. 

Let $a$ be any positive integer less than $p$ and $\sigma_a \in \mathfrak G(\mathbb Q(\zeta_{p})/\mathbb Q)$
 be the element of the Galois group of $\mathbb Q(\zeta_{p})$ whose action on $\zeta_{p}$ is defined as
  $\sigma_a(\zeta_{p})=\zeta_{p}^a$, \cite{wash},  then
 $$ \sigma(G_j) = \sum_{\chi_3( x)=\zeta_3^j} \zeta_p^{ax} =\sum_{\chi_3( x'a^{-1})=\zeta_3^j} \zeta_p^{x'}=\sum_{\chi_3( x')=\zeta_3^j\chi_3(a)} \zeta_p^{x'}  ~~. $$
 This implies that any of these automorphisms induces a permutation of $G_0$, $G_1$, and $G_2$, and therefore leaves
 their symmetric functions invariant, which thus belong to $\mathbb Q$. In particular, the three
 elementary symmetric functions
 $$  s_1=G_0+G_1+G_2 ~~,~~  s_2=G_0 G_1+ G_1 G_2 + G_2 G_1~~,~~ s_3=G_0 G_1 G_2 ~~,~~ $$
are rational integers; it follows that $G_0$, $G_1$, and $G_2$ are the roots of
 a cubic polynomial with rational coefficients
  $q(z)= z^3-s_1 z^2+s_2 z-s_3 $, which has a cyclic Galois group of order $3$
 (since $\frac{p-1}{3}$ values of $a$ give the same permutation of its roots). Thus $q(z)$ is
 irreducible over $\mathbb Q$, and denoting one root with $\eta$, the other roots can be expressed as
 polynomials with integer coefficients $r_1(\eta)$ and $r_2(\eta)$ of degree $2$ in $\eta$. 
\QED

\noindent
Gauss computed the coefficients of the cubic polynomial $q(z)$ by a clever manipulation of the periods,
 a task that generally has non-polynomial-time complexity in $p$. 
The following theorem proves that these coefficients can be computed, with deterministic
 polynomial-time complexity, exploiting pure arithmetic features 
 of $p$ without dealing with Gauss periods.

\begin{theorem}
   \label{theo3valc} 
Let $q(z)$ be the monic polynomial whose roots are $G_0$, $G_1$ and $G_2$. Then 
\begin{equation}
  \label{gausscubic}
q(z)=z^3+z^2-\frac{p-1}{3}z-\frac{(3+u)p-1}{27},
\end{equation}
 where $u$ is obtained from the
 representation $4p=u^2+27v^2$ and taken with the sign making the constant term
 an integer. 
\end{theorem} 

\noindent
{\sc Proof}.
Let $q(z)$ be as above $z^3-s_1 z^2+s_2 z-s_3 $. It is immediately seen that $s_1=-1$, as $\sum_{x=0}^{p-1} \zeta_p^{x}=0$. %, and it is natural to ask whether $s_2$ and $s_3$ can also be
 %easily computed given $p$. 
Using the structure constants of the integral algebra generated by $G_0$, $G_1$, and $G_2$, we will show that $s_2=-\frac{p-1}{3}$,
 while $s_3$ ultimately depends on the representation $(u,v)$ of $4p$ by the quadratic form
 $u^2+27v^2$.  \\
Let $\sigma$ be a generator of the cyclic Galois group of $\mathbb Q(\eta)$ ($\sigma$ is
 also a generator of the Galois group of $q(z)$), then $G_0=\eta$, $G_1=\sigma(\eta)$,
 and $G_2=\sigma^2(\eta)$ are $\mathbb Q$-linearly independent \cite{artin} and generate an
 algebra, \cite[Lemma 2.2]{monico}, whose constants of multiplication \cite{dickson} are
 rational integers, \cite[Remark 2.3]{monico}, so that 
$G_0 G_1$, $G_1 G_2$, and $G_2 G_0$ are linear combinations of $G_0$, $G_1$, and $G_2$ with integer
 coefficients. Furthermore, since $G_0$, $G_1$, and $G_2$ are cyclically permuted by the action of $\sigma$,
 we can write 
\begin{equation}
   \label{alg1}
   \left\{ \begin{array}{l}
       G_0 G_1 = a G_0 + b G_1+ c G_2  \\ 
       G_1 G_2 = a G_1 + b G_2+ c G_0  \\ 
       G_2 G_0 = a G_2 + b G_0+ c G_1  \\ 
    \end{array} \right.
\end{equation}
where $a, b, c$ are integers whose sum is $\frac{p-1}{3}$, since each $G_j$ contains $\frac{p-1}{3}$ powers of $\zeta_p$ and each $G_iG_j$, $i \neq j$, expands into $\frac{(p-1)^2}{9}$ terms which are powers of $\zeta_p$ whose exponents, reduced modulo $p$, are never $0$.
Then, summing the three equations, we get
$$s_2 = (a+b+c)(G_0+G_1+G_2) = -\frac{p-1}{3}~~.$$  
The evaluation of $s_3$ requires the explicit knowledge of $c$:
 by summing the three equations, multiplied by $G_2$, $G_0$, and $G_1$ respectively, we obtain the relation 
$$ 3 G_0 G_1 G_2 = (a+b) s_2 + c (G_0^2+G_1^2+G_2^2) = (\frac{p-1}{3}-c) s_2 + c (1-2s_2)~~,  $$
which yields $s_3= \frac{1}{3}[cp-\frac{(p-1)^2}{9}]$. Now, the value of $c$ is specified as follows. \\
%
%\begin{theorem}
 %  \label{theo3valc} 
%Let $p(z)=z^3+z^2-\frac{p-1}{3}z-\frac{1}{3}[cp-\frac{(p-1)^2}{9}]$  be the polynomial of the
% periods $G_0$, $G_1$, and $G_2$, then $c=\frac{p+1+u}{9}$, where $u$ is obtained from the
% representation $4p=u^2+27v^2$ and taken with the sign making $c$ an integer, that is,
% $u$ should be taken with the sign that corresponds to the $u$ determination modulo $9$ shown 
%in the following table 
%$$ \begin{array}{|c|c|} \hline
%     p \bmod 9 & u \bmod 9  \\  \hline
 %    1         &      7  \\
 %    7         &      1  \\
  %   4         &      4  \\  \hline
  %   \end{array}
%$$     
%\end{theorem} 
%
%\noindent
%{\sc Proof}.
Since the Galois group of $q(z)$, which is a polynomial with integer coefficients, is cyclic of order $3$,
 its discriminant $\Delta$ is the square of an integer \cite[Proposition 7.4.2]{cox}. The direct
 computation yields
 $$ \Delta= - \frac{p^2 (p^2-2p+1-18cp-18c+81c^2)}{27}  ~~, $$
then $\Delta$ must be of the form $v^2 p^2$, whence $c$ is obtained from the equation
\begin{equation}
\label{4p}
 0 = p^2-2p+1-18cp-18c+81c^2+27 u^2 =-4p+(9c-p-1)^2+27v^2 ~~. 
\end{equation}
It is known \cite{huat,gauss,hardy,rosen} that primes of the form $6k+1$ are essentially (up to signs) represented in a unique way
 by the quadratic form $x^2+3y^2$ (note that this representation may be computed  
 in deterministic polynomial time using the Schoof algorithm \cite{schoof} and the Gauss
 reduction algorithm of quadratic forms \cite{mathews}).
Further, $4=1+3$ is represented by the same form, then 
 $4p$ has essentially three different representations, namely
$$  4p=(2x)^2+3(2y)^2 ~~~~,~~~~ 4p=(x-3y)^2+3(x+y)^2     ~~~~,~~~~  4p=(x+3y)^2+3(x-y)^2      ~~~~. $$
  Since $x$ is relatively prime with $3$, necessarily exactly one of $2y$, or $x+y$ or $x-y$ is divisible
  by $3$ and allows us to write $4p=u^2+27 v^2$.  By comparison with (\ref{4p}) we have $9c-p-1 = u$,
  and $u$ should be taken with the sign that makes the expression $p+1+u$ divisible by $9$, in order
  to have an integer $c$, that is $u$ should be taken with the sign that makes it congruent to $1$
  modulo $3$: since $p+1+u=0 \bmod 9$
  implies that the same expression is $0$ modulo $3$ and $p$ is
  congruent to $1$ modulo $3$,  then also $u$ must be congruent $1$ modulo $3$. \hfill $\Box$
%$$ \begin{array}{|c|c|} \hline
%     p \bmod 9 & u \bmod 9  \\  \hline
%     1         &      7  \\
%     7         &      1  \\
%     4         &      4  \\  \hline
%     \end{array}
%$$     

%  
%\QED

\begin{corollary}
   \label{corolabc} 
The multiplicative constants of the algebra generated by the periods (equation (\ref{alg1}))
 are obtained from the representation $4p=u^2+27 v^2$ as
$$ a= \frac{2p-u+9v-4}{18} ~~~~,~~~~ b=\frac{2p-u-9v-4}{18} ~~~~,~~~~c=\frac{p+1+u}{9} , $$  
 where the sign of $u$ is specified in Theorem \ref{theo3valc}, and the sign
 of $v$ is such that $a$ and $b$ are compliant with the chosen definitions of
  $G_1$ and $G_2$.
 %since the naming of the periods is arbitrary even once $G_0$ is given.
\end{corollary} 
\noindent
{\sc Proof}.
The value $c=\frac{p+1+u}{9}$ has been found in the proof of Theorem \ref{theo3valc}, where it was also  remarked that $a+b+c=\frac{p-1}{3}$, thus for computing $a$ and $b$, 
%besides $a+b=\frac{2p-4-u}{9}$,
 we only need a further independent relation. \\
Adding member by member the three equations in (\ref{alg1}) after their orderly multiplication
 by $G_0$, $G_1$, $G_2$, or by $G_1$, $G_2$, $G_0$, respectively, we obtain 
$$ r_1=G_0^2G_1+G_1^2G_2+G_2^2G_0= a(G_0^2+G_1^2+G_2^2)+(b+c)(G_0G_1+G_1G_2+G_2G_0) ~~, $$
$$ r_2=G_1^2G_0+G_0^2G_2+G_2^2G_1= b(G_0^2+G_1^2+G_2^2)+(a+c)(G_0G_1+G_1G_2+G_2G_0) ~~. $$
% 
%The second member in both expressions is a symmetric function (for the symmetric group) of $G_0$, $G_1$,
% and $G_2$, thus it is invariant for the exchange of the index $0$ and $1$, while, the first members are 
% exchanged one each other, this exchange forces an exchange of $a$ with $b$.
If we know either $r_1$, or $r_2$, then we have a second linear equation for $a$ and $b$.
To compute $r_1$ and $r_2$, we observe that they are exchanged by permuting, for example, $G_0$ and $G_1$,
 and are invariant under a cyclic permutation of $G_0$, $G_1$, and $G_2$.
Then, their sum $r_1+r_2$ and product $r_1r_2$ are symmetric functions of the roots of $q(z)$
  and a theorem of Lagrange's \cite{burnside} assures that they can be expressed by means of the
  elementary symmetric functions $s_1$, $s_2$, and $s_3$ (the coefficients of $q(z)$).
Thus, using properties of the symmetric functions \cite{burnside}, we have
$$ r_1+r_2 = s_1 s_2-3 s_3 ~~~~,~~~~r_1r_2 = s_2^3-6 s_1s_2 s_3-9 s_3^2 -s_3 s_1^3 ~~,$$
   and substituting the explicit values of $s_1$, $s_2$, and $s_3$ given in Theorem \ref{theo3valc},
   we obtain
$$ r_1+r_2 = \frac{3p-1-p(u+3)}{9}~~~~\mbox{and}~~~~r_1r_2 = \frac{1+pu+p^2u^2-3p^3}{81} ~~. $$
Solving a second degree equation, we find $r_1= \frac{1}{2}(\frac{3p-1-p(u+3)}{9}+pv)$ and
 $r_2= \frac{1}{2}(\frac{3p-1-p(u+3)}{9}-pv)$, where the sign of $v$ should be properly chosen to match
 the values of $r_1$ and $r_2$ obtained from the definition of the Gauss periods.
In conclusion, from the system
$$  \left\{ \begin{array}{lcl}
       a+b &=&\sty \frac{2p-4-u}{9} \\
       a  \frac{2p+1}{3} -b \frac{p-1}{3} &=&\sty \frac{2p^2+27pv-pu-2u-8}{54} ~~
       \end{array} \right.  $$
 we obtain $a= \frac{2p-u+9v-4}{18}$ and $b=\frac{2p-u-9v-4}{18}$.
\QED
%%%%%%%%%%%%%%%%%%%%%% 
%\paragraph{Remark 3.}
It is possible to obtain a representation of the Gauss sum $G_{1}(1, \chi_3)$ in terms of
 a single root  $\eta_p$ of $q(z)$ by expressing $G_0$, $G_1$, and $G_2$ in terms of $\eta_p$, because $\mathbb Q(\eta_p)$ is the splitting field of $q(z)$, as seen in Theorem 3 (cf. also  \cite{artin}).
 For example, we may set $G_0=\eta_p$, thus the other roots, that is, Gauss periods, are
$$  G_1 = \frac{G_0^2}{v}+\frac{4-u-3v}{6v} G_0 +\frac{2-u-9v-4p}{18v}  ~~, $$
$$  G_2 =-\frac{G_0^2}{v}-\frac{3v-u+4}{6v} G_0 -\frac{9v-u-4p+2}{18v} ~~. $$
We note that changing the sign of $v$ is equivalent to exchanging the values of $G_1$ and $G_2$.
%
%We may remove the boring choice of the sign of $v$, or equivalently,
% the uncertainty in the attribution of the values to $a$ and $b$, that is not settled in this corollary,
% arguing as follows. \\
%Since the polynomial $q(z)$ has a cyclic Galois group, which is in particular an abelian group,
% by a theorem of Abel's              {cox},
% its roots can be expressed in terms of a single root; in other terms, the splitting field
% of $q(z)$ over $\mathbb Q$ has degree $3$ like the order of its automorphism group which is
% also the Galois group of $q(z)$ \cite{artin}.
%Therefore, the roots of $q(z)$ are quadratic polynomials in one of them, say $G_0$, 
% in particular, the other roots $G_1$ and $G_2$, are
%$$  G_1 = -\frac{G_0^2}{v}-\frac{3v-u+4}{6v} G_0 -\frac{9v-u-4p+2}{18v} ~~, $$
%$$  G_2 = \frac{G_0^2}{v}+\frac{4-u-3v}{6v} G_0 +\frac{2-u-9v-4p}{18v} ~~. $$
%Note that $\{ 1, G_0, G_0^2 \}$ is not necessarily an integral basis for the splitting field of $q(z)$,
%and this explain the denominator in the coeffcient, denominator that usually remains. \\
These equations establish a correspondence between $G_0$, and $G_1$ and $G_2$, thus 
 the coefficients $a$, $b$, and $c$ in equation (\ref{alg1}) can be uniquely specified, for
 instance the equation $G_0G_1=aG_0+bG_1+cG_2$ is satisfied choosing
$$ a= \frac{2p-u-9v-4}{18} ~~~~\mbox{and}~~~~ b=\frac{2p-u+9v-4}{18} ~~, $$
 whatever be the sign of $v$; $c$ is specified in any case as shown in Theorem \ref{theo3valc}. 

These observations yield the following representation:
\begin{theorem}
The Gauss sum $G_{1}(1, \chi_3)$
 is uniquely characterized in terms of a root $\eta_p$ of $q(z)$, with $u,v$
 obtained from the representation $u^2+27v^2$ of $4p$, as
\begin{equation} 
  \label{gauss3}
\eta_p +\zeta_3 (\frac{\eta_p^2}{v}+\frac{4-u-3v}{6v} \eta_p +\frac{2-u-9v-4p}{18v}) + \zeta_3^2 (-\frac{\eta_p^2}{v}-\frac{3v-u+4}{6v} \eta_p -\frac{9v-u-4p+2}{18v})     ~.
\end{equation}
\end{theorem}
%%%%%%%%%%%%%%%%%%%%%% 

\paragraph{Remark 3.}
Since $\zeta_3^2=-1-\zeta_3$, we can write
$$ G_{1}(1, \chi_3)= G_0- G_2 +\zeta_3 (G_1 - G_2) ~~,$$
 thus the relation $G_{1}(1, \chi_3) \bar G_{1}(1, \chi_3)=p$ yields
$$  p =(G_0- G_2 )^2-(G_0- G_2 )(G_1- G_2 )+(G_1- G_2 )^2  $$ 
which shows that the equation $x^2-xy+y^2=p$ has further solutions in the maximal order of
 $\mathbb Q(\zeta_p)$ besides the solutions in rational integers, for example $x=r_1(\eta)-\eta$
 and $y=r_2(\eta)-\eta$. 
 
\paragraph{Example}
Consider $p=7$, then the Gauss sum has the form
$$ G_1(1,\chi)= (\zeta_7+\zeta_7^6)+\zeta_3 (\zeta_7^2+\zeta_7^5)+ \zeta_3^2 (\zeta_7^3+\zeta_7^4) ~~, $$
the coefficients $G_j$ of the powers of $\zeta_3$ are real, and are roots of the cubic polynomial
$z^3+z^2-2z-1$, which has a cyclic Galois group of order $3$ over $\mathbb Q$. Let $\eta_7$ be a root of this polynomial.
 The other roots are $-2+\eta_7^2$ and $1-\eta_7-\eta_7^2$, thus if we choose the roots $\eta_7=\zeta_7+\zeta_7^6$, and the other two roots equal to $\zeta_7^2+\zeta_7^5$ and $\zeta_7^3+\zeta_7^4$, respectively, we obtain the %"closed" 
expression $\eta_7+ \zeta_3(-2+\eta_7^2)+\zeta_3^2(1-\eta_7-\eta_7^2)$,
 which coincides with the expression obtained specializing (\ref{gauss3}) with $p=7$, $u=1$, and $v=1$.
Furthermore, it is direct to check that
 $x=-2+\eta_7^2-\eta_7$ and $y=1-\eta_7-\eta_7^2-\eta_7$ give a representation of $7$ through the quadratic form
 $x^2-xy+y^2$ in integers of $\mathbb Q(\eta_7)$, which may be of interest besides the $12$ 
representations in rational integers (see e.g. \cite[Proposition 8.3.1]{rosen}, \cite{mollin}), namely $x= 2$ and $y=-1$ and those obtained through associates and conjugates of $2-\zeta_3$ in $\mathbb Q(\zeta_3)$.%its associated or conjugated. 
 %%$x=3$ and $y=-1$,  or $x=3$ and $y=-2$.

\vspace{5mm}

\noindent
%  
%  but we won't investigate about its precise value.
%.....proof....
 %$\mathbb F_{p}$ a Gauss sum of a nontrivial cubic character is an algebraic integer in
% $\mathbb Q(\zeta_{3p})$ of absolute value $\sqrt p$. 
%
A Gauss sum over $\mathbb F_{p^r}$ is in general also not rational, as can be found using
 Davenport-Hasse's theorem, \cite{berndt,Jungnickel}, by lifting the case over $\mathbb F_{p}$.
%Elementary proofs for every $r>1$ appear to require special tricks that depend on $r$ itself. 
%However, an interesting expression is readily obtained by means of
%For $r|(p-1)$ 
If there exists an irreducible polynomial $X^r-\beta$ over $\mathbb F_{p}$
we can use Lemma \ref{lemma1pbis} %considering the extension for , as proved in 
to obtain the following theorem:

\begin{theorem}
   \label{theo4bis} 
If $p=6k+1$, and if there exists an irreducible polynomial $X^r-\beta$ over $\mathbb F_{p}$, then $G_{r}(1, \chi_3)$ is again an element of the subfield
 $\mathbb Q(\zeta_{3}, \eta)$ of degree $6$ of $\mathbb Q(\zeta_{3p})$, and in particular
 it can be written in the form
$$ G_{r}(1, \chi_3) = \bar\chi_3(r) G_{1}(1, \chi_3) B_{r,1}(\alpha) ~~,  $$ 
 where $\alpha$ is a root of $X^r-\beta$. % and $\beta \in \mathbb F_p$ is a cube such that
  %$\chi_r(\beta)=\zeta_r$. 
The factor $B_{r,1}(\alpha)$ has the form $B_0 + B_1 \zeta_3 + B_2 \zeta_3^2$
  where $B_0$, $B_1$, and $B_2$ are positive rational integers, that can be computed, up to a permutation, from the solutions
  of a quadratic Diophantine equation.
\end{theorem} 

\noindent
{\sc Proof}. 
By Lemma \ref{lemma1pbis}, we have %can find such an $\alpha$, and have
$$  G_{r}(1, \chi_3)= \bar\chi_3(r)  G_{1}(1, \chi_3) B_{r,1}(\alpha)~~, $$
where $B_{r,1}(\alpha)$ is an element of $\mathbb Q(\zeta_3)$ with absolute value $\sqrt{p^{r-1}}$ (by taking absolute values of both sides). % and $\chi_3(\alpha)=1$. 
%
%\noindent 
This expression shows that $G_{r}(1, \chi_3)$ belongs to $\mathbb Q(\eta, \zeta_3)$
 as $G_{1}(1, \chi_3)$, since $B_{r,1}(\alpha)$ belongs to $\mathbb Q(\zeta_3)$. 
% $G_{r}(1, \chi_3)= \chi_3(r) G_{1}(1, \chi_3) B_r(\alpha)$ shows that
% 
The factor $B_{r,1}(\alpha) \in \mathbb Q(\zeta_3)$ can be %computed explicitely since it can be 
written as
$$   B_{r,1}(\alpha) = B_0 + B_1 \zeta_3 + B_2 \zeta_3^2 $$
where $B_0$, $B_1$, and $B_2$ are positive integers whose sum is $p^{r-1}$.
Since the square of the norm of $B_{r,1}(\alpha)$ is $p^{r-1}$, we have the Diophantine equation
$$ p^{r-1}= B_0^2+B_1^2+B_2^2 - B_0 B_1 -B_1B_2- B_2B_0 ~~.$$
Using the relation $B_0+B_1+B_2 = p^{r-1}$, we eliminate $B_2$ and obtain a quadratic Diophantine 
 equation that can be solved for $B_0$ and $B_1$:  %, it will turn out that the solution is unique:
$$ 3 B_0^2+3B_1^2+3B_0B_1-3p^{r-1}B_0-3p^{r-1}B_1+ 3(p^{2r-1}-p^{r-1}) =0 ~~.$$
With the substitution 
$$  B_0= \frac{X+p^{r-1}}{3} ~~,~~ B_1= \frac{Y+p^{r-1}}{3} ~~,  $$
we obtain the equation
$$ X^2+XY+Y^2-3p^{r-1}  =0 ~~, $$
whose solutions can be obtained from the solution of 
$$   u^2+uv+v^2= p $$
as coming from
$$ X- \zeta_3 Y =%\pm \zeta_3^ i 
(1- \zeta_3) (u- v \zeta_3)^{r-1}  $$
by composition of quadratic forms.
%where $i=0,1,2$ and the double signs will give the $12$ solutions, among which to search the right one,
%The co
The ultimate assignment of the solutions to the $B_i$ depends on the choice of the primitive roots in the definition of the Gauss sums. %\hfill $\Box$
%however, the condition that $B_i$s must be integers leaves only three possibilities that corresponds to
%a cyclic permutation of $B_i$s. In order to specifies the right correspondence, one $B_i$ should be
% evaluated, because, apparently, there is no purely arithmetic condition.

\QED

In the following we focus on the special case $r=2$ to enlighten some properties and relations of the Gauss sums seen from different perspectives.

%deserves to be considered singularly because some related results
% have an own importance and cannot be obtained simply by specializing more general results.

\begin{theorem}
   \label{theo4} 
If $p=6k+1$, then $G_{2}(1, \chi_3)$ is again an element of the subfield $\mathbb Q(\zeta_{3}, \eta)$ of
 degree $6$ of $\mathbb Q(\zeta_{3p})$, and in particular it can be written in the form
$$ G_{2}(1, \chi_3) = G_{1}(1, \chi_3) A_1(\frac{1}{2\alpha}) ~~,  $$ 
 where $\alpha$ is a root of $x^2-\beta$ and $\beta \in \mathbb F_p$ is a cube and quadratic non-residue.
\end{theorem} 

\noindent
{\sc Proof}. 
%Let $\alpha$ be defined, similarly to Lemma \ref{lemma1p}, as a root of an irreducible polynomial $X^2-\beta$, with a suitable $\beta \in  \mathbb F_{p}$, then the Gauss sum $G_{2}(1, \chi)$
% can be written as
%$$ G_{2}(1, \chi)= \sum_{z \in  \mathbb F_{p^{2}}}   \chi(z) e^{\frac{2\pi i}{p} \Tr_{2}(z)} = 
 %  \sum_{x,y \in  \mathbb F_{p}}   \chi(x+\alpha y) e^{\frac{2\pi i}{p} \Tr_{2}(x+\alpha y)}
%   =  \sum_{x,y \in  \mathbb F_{p}}   \chi(x+\alpha y) e^{\frac{2\pi i}{p} 2x} ~~. $$  
%We split the summmation into three sums
%$$ \bar \chi(2) \sum_{y \in  \mathbb F_{p}}   \chi(2\alpha y) ~~,~~
 %\bar \chi(2) \sum_{x' \in  \mathbb F_{p}^*}   \chi(x') e^{\frac{2\pi i}{p} x'} ~~,~~
%\bar \chi(2) \sum_{x',y \in  \mathbb F_{p}^*}   \chi(x'+2\alpha y) e^{\frac{2\pi i}{p} x'} ~~.$$
%The first summation is $0$ because $\chi$ is a non-trivial character in $\mathbb F_{p}$, 
 %the second summation is \\
% $\bar \chi(2) G_{1}(1, \chi)$, and the third summation
 %gives, after the substitution $y=z x'$,
%$$ \bar \chi(2) \sum_{x',z \in  \mathbb F_{p}^*}   \chi(x'+2\alpha z x') e^{\frac{2\pi i}{p} x'}=
 %   \bar \chi(2) \sum_{x' \in  \mathbb F_{p}^*}   \chi(x') e^{\frac{2\pi i}{p} x'} 
  %    \sum_{z \in  \mathbb F_{p}^*}  \chi(1+2\alpha z) = $$
%$$ \bar \chi(2) G_{1}(1, \chi) [  \chi(2\alpha) \sum_{z \in  \mathbb F_{p}}  \chi(\frac{1}{2\alpha}+z) - 1]~~. $$
As in Theorem \ref{lemma3p}, we can find such an $\alpha$ and then use Lemma \ref{lemma1p} to deduce

%Summarizing, we have
$$  G_{2}(1, \chi_3)=  G_{1}(1, \chi_3)  \sum_{z \in  \mathbb F_{p}}  \chi_3(\frac{1}{2\alpha}+z)=
     G_{1}(1, \chi_3) A_1(\frac{1}{2\alpha})~~, $$
where $A_1(\frac{1}{2\alpha})$ is an element of $\mathbb Q(\zeta_3)$ with absolute value $\sqrt{p}$ and $\chi_3(\alpha)=1$. %Then $A_1(\frac{1}{2\alpha})$ has the form
%$$ A_1(\frac{1}{2\alpha})= A_0 + A_1 \zeta_3 + A_2 \zeta_3^2  ~~,~~ A_0, A_1, A_2 \in \mathbb Z ~~.$$ 
%Now, we show that $A_1(\beta)$,  for every $\beta \in \mathbb F_{p^{2}}$ is the same, up to a permutation of $A_0$, $A_1$ and $A_2$, as $A_1(\alpha)$. \\
%Let $\beta=x+\alpha y$ be an lement of $\mathbb F_{p^2}$ with $x,y \in \mathbb F_{p}$ and $y \neq 0$,
% we have
%$$  A_1(\beta) = \sum_{z \in  \mathbb F_{p}}  \chi(x+\alpha y +z) =
% \sum_{z' \in  \mathbb F_{p}} \chi(\alpha y +z') = 
% \chi(y) \sum_{z' \in  \mathbb F_{p}} \chi(\alpha +z'y^{-1}) = \chi(y)  A_1(\alpha) ~~.$$

\noindent In conclusion $G_{2}(1, \chi_3)=G_{1}(1, \chi_3) A_1(\frac{1}{2\alpha})$ shows that
 $G_{2}(1, \chi_3)$ belongs to $\mathbb Q(\eta, \zeta_3)$.
\QED

\paragraph{Remark 4.}
Theorem \ref{theo4} states that the Gauss sum $G_{2}(1, \chi_3)$ is the product of 
  $G_{1}(1, \chi_3)$ and $A_1(\frac{1}{2\alpha})$, a result that is slightly different from that obtained
  using Davenport-Hasse's theorem \cite{rosen}, which states that
$G_{2}(1, \chi_3')= - G_{1}(1, \chi_3)^2$, where $\chi_3'$ is defined by extending the nontrivial
 character over $\mathbb F_p$ to a character over $\mathbb F_{p^m}$, using the extension rule 
$$  \chi_3'(x) = \chi_3(N_{\mathbb F_{p^m}}(x)),$$
where $N_{\mathbb F_{p^m}}(x)= x\cdot x^p\cdots x^{p^{m-1}}$ is the norm of $x$.
In our case we have $  \chi_3'(x) = \chi_3(N_{\mathbb F_{p^2}}(x))$, and whenever $\chi_3'$ is restricted to 
 $\mathbb F_p$, we specifically have
$$  \chi_3'(x) = \chi_3(N_{\mathbb F_{p^2}}(x))= \chi_3(x^2) = \bar \chi_3(x) ~~~~\forall x \in \mathbb F_p ~~.~~ $$
Therefore, since $G_{2}(1, \bar \chi_3')=\chi_3'(-1)\bar G_{2}(1, \chi_3')=\bar G_{2}(1, \chi_3')$, the equation given by Davenport-Hasse can be read as
$$ G_{2}(1, \chi_3')= - \bar G_{1}(1, \chi_3')^2 ~~,$$
where $\chi_3'$ is a cubic character defined in $\mathbb F_{p^2}$, and $\bar G_{1}(1, \chi_3')$ is
 evaluated on the subset $\mathbb F_p$.

\noindent In the following proposition we show how this relation may also be derived elementarily. First we need a well-known lemma (see also \cite[Proposition 8.3.3]{rosen} or \cite{adhikari}), for which we present an alternative proof:
\begin{lemma}
   \label{lemma5}
If $p=6k+1$, then $  G_{1}(1, \chi_3)^3= p \sum_{x \in \mathbb F_p} \chi_3(x(x-1)) $.
\end{lemma}

\noindent
{\sc Proof}. The proof is straightforward from the computation of the cube
$$ G_{1}(1, \chi_3)^3=  \sum_{x,y,z \in \mathbb F_p} \chi_3(xyz) \zeta_p^{x+y+z} =
    \sum_{x,y,u \in \mathbb F_p} \chi_3(xy(u-x-y)) \zeta_p^{u} ~~,  $$
in which the substitution $u=x+y+z$ has been performed. The summation over $x$ can be split into
 two summations $S_1$ and $S_2$, depending on whether $u=y$ or $u \neq y$. The first summation turns out to be $0$, since
$$ S_1=\sum_{y \in \mathbb F_p} \sum_{x \in \mathbb F_p} \chi_3(xy ~~ (x)) \zeta_p^{y} =
   \sum_{y \in \mathbb F_p}\chi_3(y) \zeta_p^{y}  \sum_{x \in \mathbb F_p}\chi_3(x^2)=
\sum_{y \in \mathbb F_p}\chi_3(y) \zeta_p^{y}  \sum_{x \in \mathbb F_p}\chi_3^2(x) =0 ~~.  $$ 
The second summation, with the substitution $x=x'(u-y)$, becomes
$$ S_2= \sum_{\stackrel{y, u \in \mathbb F_p}{y\neq u} } \zeta_p^{u}  \sum_{x \in \mathbb F_p}  \chi_3(xy(u-x-y))=
    \sum_{\stackrel{y, u \in \mathbb F_p}{y\neq u} }  \chi_3(y) \zeta_p^{u}  \sum_{x' \in \mathbb F_p} \bar \chi_3(u-y)  \chi_3(x'(1-x'))  ~~.  $$
Defining $A=\sum_{x' \in \mathbb F_p}\chi_3(x'(1-x')) $, a constant that does not depend on $u$ or $y$, we
may write 
$$ S_2= A \sum_{u \in \mathbb F_p} \zeta_p^{u} \sum_{y \neq u }  \chi_3(y) \bar \chi_3(u-y) =
      A[\sum_{y  \in \mathbb F_p}  \chi_3(y) \bar \chi_3(0-y) +   \sum_{u \neq 0}  \zeta_p^{u} \sum_{y \in \mathbb F_p}  \chi_3(y) \bar \chi_3(u-y) ]  ~~.   $$ 
In conclusion, we have $ S_2= Ap $, since the first summation over $y$ is $p-1$, the second summation over $y$ is $-1$ independently of $u$, \cite{gauss1,winter}; finally, the summation over $u$ is $-1$, so that $p-1+ (-1) (-1) = p$.
\QED

\noindent
Since $G_{1}(1, \chi_3) \bar G_{1}(1, \chi_3)=p$, the above result gives
$$ G_{1}(1, \chi_3)^3=G_{1}(1, \chi_3) \bar G_{1}(1, \chi_3) A $$
which implies that $G_{1}(1, \chi_3)^2= \bar G_{1}(1, \chi_3) A $. On the other hand, Theorem \ref{theo4} 
 gives 
 $$G_{2}(1, \chi_3) = G_{1}(1, \chi_3) A_1(\frac{1}{2\alpha}) ~~, $$ 
  thus we can prove the identity 
 $  G_{2}(1, \chi_3) =- \bar G_{1}(1, \chi_3)^2  $, implied by Davenport-Hasse's theorem,
 if we can prove that %the ratio $\frac{\bar A_1(\frac{1}{2\alpha})}{A}$ is a unit
 %of $\mathbb Z(\zeta_3)$, the integral ring of $\mathbb Q(\zeta_3)$, or equivalently that 
$A=-\bar A_1(\frac{1}{2\alpha})$. It is in fact immediately seen that both $A$ and $A_1(\frac{1}{2\alpha})$
 are primes of the form $a+ b \zeta_3$ and field norm $p$ in $\mathbb Z(\zeta_3)$. Less direct is 
 the exact relation between them, which we establish in the following proposition making use of the 
 function defined as
$$  F(d,i) = \sum_{\stackrel{y \in \mathbb F_p^*}{\chi_3(y)=1}}\jacobi{g^i y+d}{p} ~~,  $$ 
where $g$ is a primitive element in $\mathbb F_p$.
 % are associated primes of norm $p$ in $\mathbb Z(\zeta_3)$.

%thus it is sufficient to prove the
% equality  $(A) = (\bar A_1(\frac{1}{2\alpha}))$ between ideals, as we will do in the following
% Lemma.
 
\begin{proposition}
   \label{lemma6}
%Assuming the same hypotheses of Theorem \ref{theo4}, the previous Lemma, and that 
 %$\alpha$ is a root of the equation $x^2-\beta$ with $\beta$ a quadratic non-residue that is
 %also a perfect cube modulo $p$, then 
$$A = -\bar A_1(\frac{1}{2\alpha}) ~~.  $$ 
\end{proposition}

\noindent
{\sc Proof}. 
We can write $A$ in the following form
\begin{equation}\label{EqA}
 A= \sum_{x \in \mathbb F_{p}} \chi_3(x(x-1)) =\sum_{z \in \mathbb F_{p}} \chi_3(z^2-\frac{1}{4})~~, 
\end{equation}
 where the last expression was obtained by making the substitution $x=z+\frac{1}{2}$. Furthermore, 
 $A_1(\frac{1}{2\alpha})$ can be written in a similar form, arguing as follows:
$$ \bar A_1(\frac{1}{2\alpha}) = \sum_{x \in \mathbb F_{p}} \bar \chi_3(x+\frac{1}{2\alpha}) =
 \sum_{x \in \mathbb F_{p}} \chi_3(x+\frac{1}{2\alpha})^2 ~~. $$
%as the conjugate of $\zeta_3$ is $\zeta_3^2$. 
Furthermore, the identity
 $\chi_3(y)=\chi_3(y)^p=\chi_3(y^p)$, which is true since $p$ is congruent to $1$ modulo $3$
 and $\chi_3$ is a multiplicative character, implies
$$ \chi_3(x+\frac{1}{2\alpha})=\chi_3(x+\frac{1}{2\alpha})^p=\chi_3(x^p+\frac{1}{(2\alpha)^p})=
   \chi_3(x-\frac{1}{(2\alpha)}) ~~, $$
as $x$ and $2$ belong to $\mathbb F_p$, $\alpha$ is a root of $x^2-\beta$ and the Frobenius
 automorphism exchanges the roots.
Then
\begin{equation} \label{EqA1}
\bar A_1(\frac{1}{2\alpha}) = \sum_{x \in \mathbb F_{p}} \chi_3(x+\frac{1}{2\alpha}) \chi_3(x+\frac{1}{2\alpha}) = \sum_{x \in \mathbb F_{p}} \chi_3(x^2-\frac{1}{4\beta})~~. 
\end{equation}
We notice now that, by definition, the value of any summation $\sum_{z \in \mathbb F_{p}} \chi_3(z^2-d)$
 can be written in the form $a_0+a_1\zeta_3+a_2 \zeta_3^2$, where $a_0$, $a_1$ and $a_2$ are the
 numbers of $z\in\mathbb F_p$ such that the value of $\chi_3(z^2-d)$ is either $1$, or $\zeta_3$,
 or $\zeta_3^2$. Therefore, writing $z^2-d=g^i y$, with $\chi_3(y)=1$,  we have
$$ a_i= \sum_{\stackrel{y \in \mathbb F_p^*}{\chi_3(y)=1}} [1+ \jacobi{g^i y+d}{p}] = \frac{p-1}{3}+F(d,i)  ~~~~i=0,1,2  ~~,~~ $$
since $[1+ \jacobi{g^i y+d}{p}]$ is equal to $0$, if $g^i y+d$ is not a square;
 it is equal to $1$ if $g^i y+d=0$; and it is equal to $2$ if $g^i y+d$ is a square. \\
Then, setting $A=b_0+b_1\zeta_3+b_2 \zeta_3^2$ and 
 $\bar A_1(\frac{1}{2\alpha}) =c_0+c_1\zeta_3+c_2 \zeta_3^2$, and using the expressions for $A$ and $\bar A_1(\frac{1}{2\alpha})$ given in (\ref{EqA}) and (\ref{EqA1}), we obtain
$$ b_i= \sum_{\stackrel{y \in \mathbb F_p^*}{\chi_3(y)=1}} [1+ \jacobi{g^i y+\frac{1}{4}}{p}] ~~\mbox{and}~~
c_i= \sum_{\stackrel{y \in \mathbb F_p^*}{\chi_3(y)=1}} [1+ \jacobi{g^i y+\frac{1}{4\beta}}{p}] ~~~
i=0,1,2  ~~.~~ $$
The numbers $c_i$ can be written as follows
$$ c_i= \sum_{\stackrel{y \in \mathbb F_p^*}{\chi_3(y)=1}} [1+ \jacobi{g^i y+\frac{1}{4\beta}}{p}] = 
     \sum_{\stackrel{y \in \mathbb F_p^*}{\chi_3(y)=1}} [1- \jacobi{\beta}{p} \jacobi{g^i y+\frac{1}{4\beta}}{p}] = 
     \sum_{\stackrel{y \in \mathbb F_p^*}{\chi_3(y)=1}} [1- \jacobi{g^i \beta y+\frac{1}{4}}{p}] $$
because $\beta$ is a quadratic non-residue; furthermore, since $\beta$ is a cube, setting $w=y \beta$, we deduce that 
$$  c_i= \sum_{\stackrel{w \in \mathbb F_p^*}{\chi_3(w)=1}} [1- \jacobi{g^i w+\frac{1}{4}}{p}] = \frac{p-1}{3} -F(\frac{1}{4},i) ~~, $$
 which only differs in sign from $b_i=\frac{p-1}{3} +F(\frac{1}{4},i)$.
The proposition follows from the fact that  
$$A=b_0+b_1\zeta_3+b_2 \zeta_3^2 = (b_0-b_2)+(b_1-b_2)\zeta_3 $$ 
and 
$$\bar A_1(\frac{1}{2\alpha})=(c_0-c_2)+(c_1-c_2)\zeta_3 = 
 (b_2-b_0)+(b_2-b_1)\zeta_3 ~~. $$
\QED

\paragraph{Remark 5.}
As has been said, Gauss sums are algebraic integers that belong to a subfield of a
 cyclotomic field, and in the above theorems we found some factorizations of Gauss sums
 into elements that may belong to different subfields. For example Theorem 4 shows that
  $G_2(1,\chi_3) \in \mathbb Q(\zeta_{3}, \eta)$ can be expressed as a product of 
  $G_1(1,\chi_3) \in \mathbb Q(\zeta_{3}, \eta)$ and $A_1(\frac{1}{2\alpha}) \in \mathbb Q(\zeta_{3})$.
The general picture of the fields involved in these factorizations is shown in the following
 figure
 
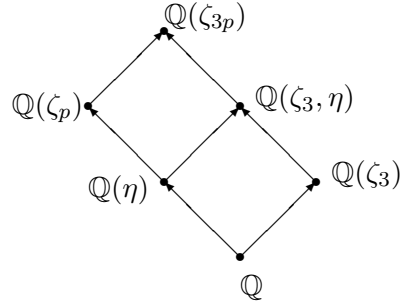
\begin{center} 
\begin{picture}(100,100)(0,0)
\setlength{\unitlength}{1mm}
\put(10,-5){$\mathbb Q$}
\put(10,0){\circle*{1}}
\put(10,0){\vector(1,1){10}}
\put(20,10){\circle*{1}}
\put(20,10){\vector(-1,1){10}}
\put(22,10){$\mathbb Q(\zeta_3)$}
\put(10,20){\circle*{1}}
\put(12,20){$\mathbb Q(\zeta_3, \eta)$}
\put(10,20){\vector(-1,1){10}}
\put( 0,30){\circle*{1}}
\put( 0,31){$\mathbb Q(\zeta_{3p})$}
\put(10,0){\vector(-1,1){10}}
\put( 0,10){\circle*{1}}
\put(-10,8){$\mathbb Q(\eta)$}
\put(0,10){\vector(1,1){10}}
\put(0,10){\vector(-1,1){10}}
\put(-10,20){\circle*{1}}
\put(-20,19){$\mathbb Q(\zeta_p)$}
\put(-10,20){\vector(1,1){10}}
%\put(70,25){\vector(1,0){10}}
%\put(80,20){\framebox(20,10){Unit 3}}
\end{picture}
\end{center}

\vspace{5mm}
\noindent
Every extension is Galois, in particular $\mathbb Q(\eta)$, $\mathbb Q(\zeta_{3})$ and $\mathbb Q(\zeta_{3},\eta)$ have Galois groups $\mathfrak G(\mathbb Q(\eta)/\mathbb Q)$,
 $\mathfrak G(\mathbb Q(\zeta_{3})/\mathbb Q)$, and  $\mathfrak G(\mathbb Q(\zeta_{3},\eta)/\mathbb Q)$,
 which are cyclic groups  of order $3$, $2$, and $6$, respectively; moreover, the third group
 $\mathfrak G(\mathbb Q(\zeta_{3},\eta)/\mathbb Q)= \mathfrak G(\mathbb Q(\zeta_{3})/\mathbb Q)\times \mathfrak G(\mathbb Q(\eta)/\mathbb Q)$ is a direct product of the other two (see also \cite{wash}). % the group of order 6 is cyclic, as it is not $S_3$, being a subgroup of an abelian group
In these fields, every rational prime $p$ of the form $6k+1$ splits into prime ideals %(\cite[pg. 138]{dedekind} or \cite[pg. 15]{wash})
 as follows:
 
\begin{center}
\begin{tabular}{l}
     $(p)=\mathfrak p^3$ in $\mathbb Q(\eta)$, i.e. the ideal $(p)$ fully ramifies; \\%p e' l'unico primo che ramifica in Q(\eta), perche' i primi che ramificano dividono il discrimante del polinomio di \eta, il quale divide il discrimante del polinomio ciclotomico di \zeta_p, visto che Q(\eta) e' contenuto in Q(\zeta_p$. Ma il discriminante del polinomio ciclotoimco e' p^{p-1}. Inoltre il polinomio di cui eta e' radice e' di terzo grado. D'altra parte \sum_i e-i f_i=m (grado del polinomio che ha generato il campo) e in questo caso abbiamo gli stessi valori per ogni i sia degli indici di ramificazione (esponenti degli ideali primi) sia dei gradi di inerzia (gradi dei polinomi nella fattorizzazone del polinomio che ha generato il campo), dato che il gruppo di galois e' ciclico (cebotarev/mollin) o semplicemente perche' gli automorfismi mandano uno nell'altro), da cui segue che e_i deve essere 3, altrimenti p sarebbe inerte e non ramificherebbe.
       \\
     $(p)=(\pi_1) (\pi_2)$ in $\mathbb Q(\zeta_3)$, i.e. the ideal $(p)$ fully splits into principal ideals; \\% class number 1, polinomio di grado 2, 3 e' l'unico prmo che ramifica e (-3||6k+1)legendre=1 quindi il polinomio splitta. -3 e' discr di x^2+x+1
       \\
     $(p)=\mathfrak P_1^3 \mathfrak P_2^3$ in $\mathbb Q(\zeta_3, \eta)$, i.e. the ideal $(p)$
       fully splits into ramified ideals; \\
       \\
     $(\pi_1)=\mathfrak P_1^3$ and  $(\pi_2)=\mathfrak P_2^3$, i.e. the principal ideals of
       $\mathbb Q(\zeta_{3})$ fully ramify in $\mathbb Q(\zeta_3, \eta)$; \\
	\\
     $\mathfrak p=\mathfrak P_1\mathfrak P_2$ in $\mathbb Q(\zeta_3, \eta)$.
\end{tabular}
\end{center}

\noindent
%The correspondence between $\pi_i$ and $\mathfrak P_i$, $i=1,2$, can be uniquely established
These factorizations can be established by the properties given in \cite[pg. 137-138]{dedekind}, that is Dedekind's formulation
 in terms of ideals of a theorem of Kummer's, or in \cite[pg. 15]{wash}. \\ 
Let $\tau_2$ denote the automorphism of order $2$ in $\mathfrak G(\mathbb Q(\zeta_{3})/\mathbb Q)$,
 which leaves the elements of $\mathbb Q(\eta)$ invariant when considered as elements of
 $\mathfrak G(\mathbb Q(\zeta_{3},\eta)/\mathbb Q)$, then $\tau_2(\mathfrak P_1)= \mathfrak P_2$.

\noindent Now, the Gauss sum $G_1(1,\chi_3)$ is an element of $\mathbb Q(\zeta_3, \eta)$ that divides $p$, as $G_1(1,\chi_3)\bar G_1(1,\chi_3)=p$, \cite{berndt}, so that $(G_1(1,\chi_3)) (\tau_2(G_1(1,\chi_3)))=(G_1(1,\chi_3)) (\bar G_1(1,\chi_3)) =(p)$.%G_1(1,\chi)=G_0+\zeta_3 G_1+\zeta_3^2 G_2
  Therefore the principal ideal $(G_1(1,\chi_3))$ will be a product of powers of the two primes $\mathfrak P_1$
 and $\mathfrak P_2$, i.e $(G_1(1,\chi_3))=\mathfrak P_1^{a}\mathfrak P_2^{b}$, where $a+b=3$ by the unique factorization in prime ideals, since the previous relation gives $(p)=\mathfrak P_1^{a}\mathfrak P_2^{b}\tau_2(\mathfrak P_1^{a}\mathfrak P_2^{b})=\mathfrak P_1^{a}\mathfrak P_2^{b}\mathfrak P_2^{a}\mathfrak P_1^{b}=\mathfrak P_1^{a+b}\mathfrak P_2^{a+b}$.

\noindent Thus, we may assume that $(G_1(1,\chi_3))=\mathfrak P_1\mathfrak P_2^{2}$, as $G_1(1,\chi_3)$ belongs properly to $\mathbb Q(\zeta_{3},\eta)$, whence Theorem \ref{theo4} and Proposition \ref{lemma6} show that $(G_2(1,\chi_3))=\mathfrak P_1^4\mathfrak P_2^{2}=(\pi_1)\mathfrak P_1\mathfrak P_2^{2}$. \\
In this framework, if the character $\chi_3'$ is used, the role of the two prime ideals is
 simply exchanged, i.e.
  $(G_2(1,\chi_3'))=\mathfrak P_2^4\mathfrak P_1^{2}=(\pi_2) \mathfrak P_2\mathfrak P_1^{2}$,
  which is the expression defined by Davenport-Hasse's theorem written in terms of 
  ideals.
  
%\vspace{2mm}
\noindent In general, the Gauss sums $G_s(1,\chi_3)$ for any $s$ can be expressed in terms of ideals as follows: $  (G_s(1,\chi_3)) = \mathfrak P_2^s\mathfrak P_1^{2s}  $.
%Moreover, the Gauss sums $G_s(1,\chi_3)$ for any odd $s$ can be expressed in terms of the ideals 
% $\mathfrak P_1$ and $\mathfrak P_2$: since $(G_s(1,\chi_3)) (\tau_2(G_s(1,\chi_3)))=(p^s)$,
% we can deduce from the Stickelberger Relation (\cite[Theorem 2, pg. 209]{rosen})
%$$  (G_s(1,\chi_3)) = \mathfrak P_1^s\mathfrak P_2^{2s} ~~. $$
 
\noindent 
However, these formulations in terms of ideals (see also \cite{denomme,katre,stickelberger}) conceal the information about which units are involved. In this sense, %the correspondence 
 %between the generators of principal ideals. 
the elementary direct approach can be more informative,
 although it may require different approaches for different situations. %special tricks for different $s$.
Considering for example the Gauss sum mentioned above, $G_1(1,\chi)$ for $p=7$ (see also \cite{garrett}),
  setting $\eta_7=\zeta_7+\zeta_7^6$, we can explicitly write the expression (which can also be obtained specializing (\ref{gauss3}) with $p=7, u=1, v=1$)
$$
G_1(1,\chi)= \eta_7+\zeta_3 (-2+\eta_7^2)+ \zeta_3^2 (1-\eta_7-\eta_7^2) ~~,
$$
whereas, choosing the ideals $\mathfrak P_1=(\zeta_3-\eta_7)$, $\mathfrak P_2=(\zeta_3^2-\eta_7)$, we must
 find a unit in order to obtain a complete factorization:
$G_1(1,\chi)=(4-\eta_7-2\eta_7^2) (\zeta_3-\eta_7) (\zeta_3^2-\eta_7)^2$, where $4-\eta_7-2\eta_7^2$ is a unit.
% questa unita' non e' un cubo e quindi non puo' essere ripartito fra i fattori

% Q(\zeta_7) e Q(\eta) hanno classe 1, quindi sono a fattorizzazione unica e gli ideali sono tutti principali
% l'ideale in Q(\eta) potrebbe essere non principale (corrisponde a una fattorizzazione in primi che non stanno nel campo) 

 %%$x=3$ and $y=-1$,  or $x=3$ and $y=-2$.

%  which above theorem can also be proved using (\cite[Theorem 5.16]{lidl}), 

%\noindent
\section*{Acknowledgment}
The Research was supported in part by the Swiss National Science
Foundation under grant No. 132256.

%*****************************************************************

\end{document}